\documentclass[11pt]{article}
\usepackage{amsmath}
\usepackage{amssymb}
\usepackage{amsthm}

\def\spd#1#2#3{{\displaystyle\frac{\partial^2 #1}{\partial #2\partial #3}}}
\def\fpd#1#2{{\displaystyle\frac{\partial #1}{\partial #2}}}

\def\dddotq{\raisebox{-1.4pt}{$\stackrel{\raisebox{-0.7pt}{.\kern-1pt.\kern-1pt.}}{q}$}}
\def\dddotqc{\raisebox{-1.4pt}{${\stackrel{\raisebox{-0.7pt}{.\kern-1pt.\kern-1pt.}}{q}}^c$}}

\def\onehalf{{\textstyle\frac12}}

\newtheorem*{myprop}{Proposition}

\begin{document}



\title{On the generalized Helmholtz conditions for Lagrangian systems with dissipative forces}

\author{M.\ Crampin, T.\ Mestdag and W.\ Sarlet\\[2mm]
{\small Department of Mathematics, Ghent University}\\
{\small Krijgslaan 281, B-9000 Ghent, Belgium}}


\date{}
\maketitle

\begin{quote}
{\small {\bf Abstract.}
In two recent papers necessary and sufficient conditions for a given
system of second-order ordinary differential equations to be of
Lagrangian form with additional dissipative forces were derived.  We
point out that these conditions are not independent and prove a stronger result
accordingly.
}
\end{quote}

\section{Introduction}

The Helmholtz conditions, for the purposes of this paper, are the
necessary and sufficient conditions for a given system of
second-order ordinary differential equations
$f_a(\ddot{q},\dot{q},q,t)=0$ to be of Euler-Lagrange type, that is,
for there to exist a Lagrangian $\Lambda(\dot{q},q,t)$ such that
\begin{equation}
f_a=\frac{d}{dt}\left(\fpd{\Lambda}{\dot{q}^a}\right)-
\fpd{\Lambda}{q^a}.\label{EL}
\end{equation}
Here $q^a$ are the generalized coordinates (collectively abbreviated
to $q$), $\dot{q}^a$ the corresponding generalized velocities, and so
on.  We shall state the conditions shortly.  The
Lagrangian is supposed to be of first order,  that is,
independent of $\ddot{q}$ and higher-order derivative (or more
properly jet) coordinates.

In two recent papers the problem of finding analogous necessary and
sufficient conditions for a given set of functions $f_a(\ddot{q},\dot{q},q,t)$ to take the
more general form
\begin{equation}
f_a=\frac{d}{dt}\left(\fpd{\Lambda}{\dot{q}^a}\right)-
\fpd{\Lambda}{q^a}+\fpd{D}{\dot{q}^a}\label{diss}
\end{equation}
for first-order functions $\Lambda$ (a Lagrangian) and $D$ (a
dissipation function) has been discussed.  We shall say that in this
case the equations $f_a=0$ are of Lagrangian form with dissipative
forces of gradient type.  A set of necessary
and sufficient conditions is given, in terms of standard coordinates,
in the fairly recent paper \cite{germ1}.  In a very recent paper
\cite{germ2} a version of the conditions expressed in terms of
quasi-velocities, or as the authors call them nonholonomic velocities,
is obtained.  We shall quote the conditions from \cite{germ1}
explicitly in Section~2.  These conditions are described as
generalized Helmholtz conditions to distinguish them from the
Helmholtz conditions discussed in our opening paragraph, which may be
called the classical Helmholtz conditions; these must of course
comprise a special case of the generalized conditions.

The main purpose of the present paper is to point out that the
generalized Helmholtz conditions as stated in \cite{germ1} are not
independent:\ in fact two of them are redundant, in that they can be
derived from the remaining ones.  This we show in Section~2 below.  We
use the same formalism as \cite{germ1}. Since the version of
the generalized conditions obtained in \cite{germ2} is equivalent to
that in \cite{germ1}  the same redundancy
is present there as well.  By taking advantage of the improvement in the formulation of the
generalized Helmholtz conditions that we achieve, we are able to give
a shorter and more elegant proof of their sufficiency than the one to
be found in \cite{germ1}.

There are in fact several interesting questions raised by the two
papers \cite{germ2,germ1}, only one of which will be dealt with here.
In the third and final section of our paper we give an outline of
these additional points of interest, which will receive a fuller
airing elsewhere.

We employ the Einstein summation convention throughout.

We end this introduction with a brief summary of the results about the
classical Helmholtz conditions that we shall need.

The classical Helmholtz conditions are that the $f_a$ should satisfy
\begin{align}
\fpd{f_a}{\ddot{q}^b}&=\fpd{f_b}{\ddot{q}^a}\label{HC1}\\
\fpd{f_a}{\dot{q}^b}+\fpd{f_b}{\dot{q}^a}
&=2\frac{d}{dt}\left(\fpd{f_b}{\ddot{q}^a}\right)\label{HC2}\\
\fpd{f_a}{q^b}-\fpd{f_b}{q^a}&=
\onehalf\frac{d}{dt}\left(\fpd{f_a}{\dot{q}^b}-\fpd{f_b}{\dot{q}^a}\right).\label{HC3}
\end{align}
It is a consequence of these (and not an extra condition, as stated in
\cite{germ1}) that
\[
\spd{f_a}{\ddot{q}^b}{\ddot{q}^c}=0
\]
(this follows from the vanishing of the coefficient of \dddotqc\
in condition~(\ref{HC2})).  Thus we may write
$f_a=g_{ab}\ddot{q}^b+h_a$, the coefficients being of first order,
with $g_{ab}=g_{ba}$ as a result of condition~(\ref{HC1}).  The
Helmholtz conditions can be re-expressed in terms of $g_{ab}$ (assumed
to be symmetric) and $h_a$, when they reduce to the following three
conditions:
\begin{align}
\fpd{g_{ab}}{\dot{q}^c}-\fpd{g_{ac}}{\dot{q}^b}&=0\label{2gh1}\\
\fpd{h_a}{\dot{q}^b}+\fpd{h_b}{\dot{q}^a}
&=2\frac{\bar{d}}{dt}(g_{ab})\label{2gh2}\\
2\left(\fpd{h_a}{q^b}-\fpd{h_b}{q^a}\right)
&=\frac{\bar{d}}{dt}\left(\fpd{h_a}{\dot{q}^b}-\fpd{h_b}{\dot{q}^a}\right),\label{2gh3}
\end{align}
where
\[
\frac{\bar{d}}{dt}=\fpd{}{t}+\dot{q}^c\fpd{}{q^c}.
\]
That is to say, the
$f_a$ are the Euler-Lagrange expressions of some first-order
Lagrangian if and only if $f_a=g_{ab}\ddot{q}^b+h_a$ for some
first-order functions $g_{ab}$ and $h_a$ such that $g_{ab}=g_{ba}$ and
(\ref{2gh1})--(\ref{2gh3}) hold.  This reformulation can be found in
the book by Santilli \cite{Sant}, for example.

\section{The generalized Helmholtz conditions}

We next turn to the analysis of
the generalized Helmholtz conditions.  Following the notation of
\cite{germ1} we set
\begin{align*}
r_{ab}&=\fpd{f_a}{q^b}-\fpd{f_b}{q^a}
+\onehalf\frac{d}{dt}\left(\fpd{f_b}{\dot{q}^a}-\fpd{f_a}{\dot{q}^b}\right)\\
s_{ab}&=\onehalf\left(\fpd{f_a}{\dot{q}^b}+\fpd{f_b}{\dot{q}^a}\right)
-\frac{d}{dt}\left(\fpd{f_b}{\ddot{q}^a}\right).
\end{align*}
The generalized Helmholtz conditions as given in \cite{germ1} are that
$r_{ab}$ and $s_{ab}$ are of first order, and that in addition
\begin{align}
\fpd{f_a}{\ddot{q}^b}&=\fpd{f_b}{\ddot{q}^a}\label{GHC1}\\
\fpd{s_{ab}}{\dot{q}^c}&=\fpd{s_{ac}}{\dot{q}^b}\label{GHC2}\\
\fpd{r_{ab}}{\dot{q}^c}&=\fpd{s_{ac}}{q^b}-\fpd{s_{bc}}{q^a}\label{GHC3}\\
0&=\fpd{r_{ab}}{q^c}+\fpd{r_{bc}}{q^a}+\fpd{r_{ca}}{q^b}.\label{GHC4}
\end{align}
Our main concern will be with analyzing conditions~(\ref{GHC2})--(\ref{GHC4}), which
correspond to (2.3e), (2.3f) and (2.3g) of \cite{germ1}; we shall show
that conditions (\ref{GHC2}) and (\ref{GHC4}) are redundant, being
consequences of the remaining conditions.

Our first aim is to understand exactly what it means for $r_{ab}$ and
$s_{ab}$ to be of first order, bearing in mind condition~(\ref{GHC1})
above.  From the vanishing of the coefficient of \dddotqc\ in $s_{ab}$ we have
\[
\spd{f_a}{\ddot{q}^b}{\ddot{q}^c}=0.
\]
As before $f_a=g_{ab}\ddot{q}^b+h_a$, the coefficients being of first order and
$g_{ab}$ symmetric.  The coefficient of \dddotqc\ in $r_{ab}$ is
\[
\onehalf\fpd{}{\ddot{q}^c}\left(
\fpd{g_{bd}}{\dot{q}^a}\ddot{q}^d+\fpd{h_b}{\dot{q}^a}
-\fpd{g_{ad}}{\dot{q}^b}\ddot{q}^d-\fpd{h_a}{\dot{q}^b}\right),
\]
whence
\[
\fpd{g_{bc}}{\dot{q}^a}=\fpd{g_{ac}}{\dot{q}^b}.
\]
The coefficient of $\ddot{q}^c$ in $s_{ab}$, namely
\[
\fpd{g_{ac}}{\dot{q}^b}+\fpd{g_{bc}}{\dot{q}^a}-2\fpd{g_{ab}}{\dot{q}^c},
\]
vanishes as a consequence.  The coefficient of $\ddot{q}^c$ in
$r_{ab}$ is
\[
\fpd{g_{ac}}{q^b}-\fpd{g_{bc}}{q^a}
-\onehalf\left(\spd{h_a}{\dot{q}^b}{\dot{q}^c}-\spd{h_b}{\dot{q}^a}{\dot{q}^c}\right),
\]
an expression which for later convenience we write as $\rho_{abc}$; we
must of course have $\rho_{abc}=0$.  The remaining terms in $r_{ab}$
and $s_{ab}$ are all of first order, and we have
\begin{align*}
r_{ab}&=\fpd{h_a}{q^b}-\fpd{h_b}{q^a}
-\onehalf\frac{\bar{d}}{dt}\left(\fpd{h_a}{\dot{q}^b}-\fpd{h_b}{\dot{q}^a}\right)\\
s_{ab}&=\onehalf\left(\fpd{h_a}{\dot{q}^b}+\fpd{h_b}{\dot{q}^a}\right)
-\frac{\bar{d}}{dt}(g_{ab});
\end{align*}
compare with (\ref{2gh2}) and (\ref{2gh3}), and also with equations
(2.16b) and (2.17c) of \cite{germ1}.

The redundancy of condition~(\ref{GHC2}) is a consequence of the
vanishing of $\rho_{abc}$, as we now show.  We have
\[
\fpd{s_{ac}}{\dot{q}^b}
=\onehalf\left(\spd{h_a}{\dot{q}^b}{\dot{q}^c}+\spd{h_c}{\dot{q}^a}{\dot{q}^b}\right)
-\frac{\bar{d}}{dt}\left(\fpd{g_{ac}}{\dot{q}^b}\right)-\fpd{g_{ac}}{q^b},
\]
using the commutation relation
\begin{equation}
\left[\fpd{}{{\dot q}^a} , \frac{\bar{d}}{dt} \right] =
\fpd{}{q^a}.\label{comm}
\end{equation}
It follows that $\partial s_{ac}/\partial \dot{q}^b-\partial s_{bc}/\partial \dot{q}^a=-\rho_{abc}=0$.
That is to say, condition~(\ref{GHC2}) holds automatically as a
consequence of the first-order property.  Furthermore, $\rho_{abc}=0$ is equivalent
to equation (2.17b) of \cite{germ1}; in other words, the redundancy of
(\ref{GHC2}) is actually implicit in \cite{germ1}, though not
apparently recognized there.

Before proceeding to consider condition~(\ref{GHC4}) we turn aside to
make some remarks about the classical Helmholtz conditions. The calculations just carried out are essentially the same as those
which lead to the version of the classical Helmholtz conditions given
in equations~(\ref{2gh1})--(\ref{2gh3}) at the end of the
introduction.  It is easy to see that in that case $r_{ab}=s_{ab}=0$
are necessary conditions.  This observation, together with the part of
the argument concerning the vanishing of the coefficients of \dddotq\
and $\ddot{q}$, leads to the following conditions:
\begin{align}
\fpd{g_{ac}}{\dot{q}^b}+\fpd{g_{bc}}{\dot{q}^a}&=2\fpd{g_{ab}}{\dot{q}^c}\label{gh1}\\
\fpd{h_a}{\dot{q}^b}+\fpd{h_b}{\dot{q}^a}
&=2\frac{\bar{d}}{dt}(g_{ab})\label{gh2}\\
\fpd{g_{ab}}{\dot{q}^c}-\fpd{g_{ac}}{\dot{q}^b}&=0\label{gh3}\\
\spd{h_a}{\dot{q}^b}{\dot{q}^c}-\spd{h_b}{\dot{q}^a}{\dot{q}^c}
&=2\left(\fpd{g_{ac}}{q^b}-\fpd{g_{bc}}{q^a}\right)\label{gh4}\\
2\left(\fpd{h_a}{q^b}-\fpd{h_b}{q^a}\right)
&=\frac{\bar{d}}{dt}\left(\fpd{h_a}{\dot{q}^b}-\fpd{h_b}{\dot{q}^a}\right)\label{gh5}.
\end{align}
These are the conditions quoted in Remark~3 of Section~1 of
\cite{germ1}.  However, it is now evident that two of them are
redundant.  Clearly condition~(\ref{gh1}) (which is the vanishing of
the coefficient of $\ddot{q}^c$ in $s_{ab}$) follows from
condition~(\ref{gh3}) and the symmetry of $g_{ab}$.
Condition~(\ref{gh4}) is the condition $\rho_{abc}=0$.  The second
part of the argument above, that leading to the relation
$\rho_{abc}=\partial s_{bc}/\partial\dot{q}^a-\partial
s_{ac}/\partial\dot{q}^b$, shows that in the classical case
condition~(\ref{gh4}) follows from the other conditions.  When
these two redundant conditions are removed we obtain the classical
Helmholtz conditions in the form given at the end of the introduction.

These results in the classical case are actually very well known,
though not apparently to the authors of \cite{germ1}, and have been
known for a long time:\ they are to be found, for example, in
Santilli's book of 1978 \cite{Sant} (which is in fact referred to in
\cite{germ1}).
For the sake of clarity we should point out a difference
between the two cases:\ in the classical case condition~(\ref{gh4}) is
completely redundant; in the generalized case it is not redundant, but
occurs twice in the formulation of the conditions in \cite{germ1},
once in the requirement that $r_{ab}$ should be of first order and
once as the condition $\partial s_{ab}/\partial\dot{q}^c=\partial
s_{bc}/\partial\dot{q}^a$.

We now return to the generalized conditions, and prove that
condition~(\ref{GHC4}) follows from condition~(\ref{GHC3}).  It will
be convenient to write condition~(\ref{GHC4}) as
\[
\sum_{a,b,c} \fpd{r_{ab}}{q^c}=0,
\]
where $\sum_{a,b,c}$ stands for the cyclic sum over $a$, $b$ and $c$, here
and below. As a preliminary remark,
note that if $k_{abc}$ is symmetric in $b$ and $c$ (say) then
$\sum_{a,b,c}k_{abc}=\sum_{a,b,c}k_{bac}$. Now
\[
\fpd{r_{ab}}{q^c}=\spd{h_a}{q^b}{q^c}-\spd{h_b}{q^a}{q^c}
-\onehalf\frac{\bar{d}}{dt}\left(\spd{h_a}{q^c}{\dot{q}^b}-\spd{h_b}{q^c}{\dot{q}^a}\right),
\]
and so by the preliminary remark
\[
\sum_{a,b,c}\fpd{r_{ab}}{q^c}=
-\onehalf\frac{\bar{d}}{dt}\left(\sum_{a,b,c}
\left(\spd{h_a}{q^c}{\dot{q}^b}-\spd{h_b}{q^c}{\dot{q}^a}\right)\right)
=
-\onehalf\frac{\bar{d}}{dt}\left(\sum_{a,b,c}
\left(\spd{h_a}{q^c}{\dot{q}^b}-\spd{h_a}{q^b}{\dot{q}^c}\right)\right).
\]
On the other hand, using the commutation relation (\ref{comm}) and the fact
that $\partial/\partial q^a$ and $\bar{d}/dt$ commute it is easy
to see that condition~(\ref{GHC3}) leads to
\[
\onehalf\sum_{a,b,c}\left(\spd{h_a}{q^c}{\dot{q}^b}-\spd{h_a}{q^b}{\dot{q}^c}\right)
=-\frac{\bar{d}}{dt}(\rho_{abc})=0.
\]

We therefore reach the following proposition, which is
stronger than the corresponding result in
\cite{germ1}.

\begin{myprop} The necessary and sufficient conditions for the
equations $f_a(\ddot{q},\dot{q},q,t)=0$ to be of Lagrangian form with
dissipative forces of gradient type as in (\ref{diss}) are that the
functions $r_{ab}$ and $s_{ab}$ are of first order, that
\[
\fpd{f_a}{\ddot{q}^b}=\fpd{f_b}{\ddot{q}^a},
\]
and that
\begin{equation}
\fpd{r_{ab}}{\dot{q}^c}=\fpd{s_{ac}}{q^b}-\fpd{s_{bc}}{q^a}\label{prop}.
\end{equation}
\end{myprop}

Just as in the classical case we can give an equivalent formulation of
these conditions in terms of $g_{ab}$ and $h_a$.  Bearing in mind that
$r_{ab}$ and $s_{ab}$ being of first order are essential hypotheses,
we find that the following conditions are equivalent to those given in
the proposition above:\ $f_a=g_{ab}\ddot{q}^b+h_a$ with $g_{ab}$
symmetric, where $g_{ab}$, $h_a$ are of first order and further satisfy
\begin{align*}
\fpd{g_{ab}}{\dot{q}^c}-\fpd{g_{ac}}{\dot{q}^b}&=0\\
\spd{h_a}{\dot{q}^b}{\dot{q}^c}-\spd{h_b}{\dot{q}^a}{\dot{q}^c}
&=2\left(\fpd{g_{ac}}{q^b}-\fpd{g_{bc}}{q^a}\right)\\
\sum_{a,b,c}
\left(\spd{h_a}{q^b}{\dot{q}^c}-\spd{h_a}{q^c}{\dot{q}^b}\right)&=0.
\end{align*}
The first of these is one of the classical conditions.  The second is
the condition $\rho_{abc}=0$, which holds in the classical case as we
have shown.  The third is just condition~(\ref{prop}) above expressed
in terms of $g_{ab}$ and $h_a$ (or as it turns out, in terms of $h_a$
alone), and $r_{ab}=s_{ab}=0$ in the classical case.  It is evident
therefore that the conditions above are indeed a generalization of
those for the classical case.

We end this section by giving an
alternative proof of the sufficiency of the generalized Helmholtz
conditions, based on this formulation of them, which is shorter and in
our view more elegant than the proof in \cite{germ1} (necessity is an
easy if tedious calculation).

We note first that if $g_{ab}$ is symmetric and satisfies
$\partial g_{ab}/\partial\dot{q}^c=\partial g_{ac}/\partial\dot{q}^b$
then
\[
g_{ab}=\spd{K}{\dot{q}^a}{\dot{q}^b}
\]
for some function $K=K(\dot{q},q,t)$ (a well-known result, which also
appears in \cite{germ1}).  Of course $K$ is not determined by this
relation; in fact if $\Lambda=K+P_a\dot{q}^a+Q$, where $P_a$ and $Q$
are any functions of $q$ and $t$, then $\Lambda$ has the same Hessian
as $K$ (the same $g_{ab}$, in other words).  Our aim is to choose
$P_a$ and $Q$ so that the given equations are of Lagrangian form with
dissipative forces of gradient type as in (\ref{diss}), with
Lagrangian $\Lambda$, assuming that the generalized Helmholtz
conditions above hold.  In fact we won't need to consider $Q$ because
it can be absorbed:\ if $\Lambda$ is a Lagrangian and $D$ a
dissipation function for some functions $f_a$, so are $\Lambda+Q$ and
$D+\dot{q}^a\partial Q/\partial q^a$.  We shall therefore take $Q=0$
below.

Let $E_a$ be the Euler-Lagrange expressions of $K$.
Then $E_a=g_{ab}\ddot{q}^b+k_a$ for some first-order $k_a$, by
construction, so $f_a-E_a=h_a-k_a=\kappa_a$ say, where $\kappa_a$ is
also of first order.  Moreover, $f_a$ satisfies the generalized
Helmholtz conditions by assumption, and $E_a$ does so by construction
(it satisfies the classical conditions after all), whence $\kappa_a$
satisfies
\begin{align}
\spd{\kappa_a}{\dot{q}^b}{\dot{q}^c}-\spd{\kappa_b}{\dot{q}^a}{\dot{q}^c}
&=0\label{gkappa1}\\
\sum_{a,b,c}
\left(\spd{\kappa_a}{q^b}{\dot{q}^c}-\spd{\kappa_a}{q^c}{\dot{q}^b}\right)&=0.\label{gkappa2}
\end{align}
Let us set $\partial\kappa_a/\partial\dot{q}^b-\partial\kappa_b/\partial\dot{q}^a=R_{ab}$.
Then by (\ref{gkappa1}) $R_{ab}$ is independent of $\dot{q}$, and by
(\ref{gkappa2})
\[
\sum_{a,b,c} \fpd{R_{ab}}{q^c}=0.
\]
There are therefore functions $P_a(q,t)$ such that
\[
R_{ab}=2\left(\fpd{P_a}{q^b}-\fpd{P_b}{q^a}\right)=
\fpd{\kappa_a}{\dot{q}^b}-\fpd{\kappa_b}{\dot{q}^a},
\]
which is to say that if we set
\[
\pi_{ab}=
\fpd{\kappa_a}{\dot{q}^b}-\left(\fpd{P_a}{q^b}-\fpd{P_b}{q^a}\right)
\]
then $\pi_{ba}=\pi_{ab}$. Moreover,
\[
\fpd{\pi_{ab}}{\dot{q}^c}=\spd{\kappa_a}{\dot{q}^b}{\dot{q}^c}=\fpd{\pi_{ac}}{\dot{q}^b}.
\]
It follows (just as is the case for $g_{ab}$) that there is a
first-order function $D'$ such that
\[
\pi_{ab}=\spd{D'}{\dot{q}^a}{\dot{q}^b},
\]
from which we obtain
\[
\kappa_a=\left(\fpd{P_a}{q^b}-\fpd{P_b}{q^a}\right)\dot{q}^b+\fpd{D'}{\dot{q}^a}+S_a
\]
where $S_a$ is independent of $\dot{q}$. Now take $\Lambda=K+P_a\dot{q}^a$.
Denoting the Euler-Lagrange expressions of $K$ by $E_a$ as before, the
Euler-Lagrange expressions for $\Lambda$ are
\[
E_a+\left(\fpd{P_a}{q^b}-\fpd{P_b}{q^a}\right)\dot{q}^b+\fpd{P_a}{t}
=E_a+\kappa_a-\fpd{D'}{\dot{q}^a}-S_a+\fpd{P_a}{t}.
\]
Thus, putting
\[
D=D'+\left(S_a-\fpd{P_a}{t}\right)\dot{q}^a,
\]
we get
\[
f_a=\frac{d}{dt}\left(\fpd{\Lambda}{\dot{q}^a}\right)-
\fpd{\Lambda}{q^a}+\fpd{D}{\dot{q}^a}
\]
as required.

This method of proof works equally well in the classical
case.  The proof is constructive, in the
same sense that the one in \cite{germ1} is, in either case.  It is
particularly well adapted to the familiar situation in which $g_{ab}$
is independent of $\dot{q}$, when one can take the kinetic energy
$\onehalf g_{ab}\dot{q}^a\dot{q}^b$ for $K$.

\section{Concluding remarks}

We wish to make four remarks in conclusion.

The first remark concerns the nature of
conditions~(\ref{GHC2})--(\ref{GHC4}) on the derivatives of $r_{ab}$
and $s_{ab}$, as originally expressed in \cite{germ1} (that is,
ignoring the question of dependence).  In particular, bearing in mind
the fact that $r_{ab}$ is skew in its indices, the condition
\[
\sum_{a,b,c} \fpd{r_{ab}}{q^c}=0
\]
is suggestive:\ if perchance the $r_{ab}$ were functions of the $q$ alone this
would have a natural interpretation in terms of the exterior calculus,
being the condition for the 2-form $r_{ab}dq^a\wedge dq^b$ to be
closed, that is, to satisfy $d(r_{ab}dq^a\wedge dq^b)=0$.  This point
is made, in somewhat different terms, in \cite{germ1} (and we appealed
to the same general result in our proof of sufficiency of the
generalized Helmholtz conditions in Section 2).  The authors of
\cite{germ1} go on to say, however, that the condition above `can be
interpreted as the vanishing curvature of a symplectic space', which
seems to us not to be entirely convincing.  In fact it is possible to
interpret the three conditions~(\ref{GHC2})--(\ref{GHC4}) collectively
as signifying the vanishing of a certain exterior derivative of a
certain 2-form on the space of coordinates
$t,q,\dot{q},\ddot{q},\ldots$\,, a 2-form whose coefficients involve
both $r_{ab}$ and $s_{ab}$.  This interpretation really arises from
seeing the problem in the context of the so-called variational
bicomplex (see \cite{vitolo} for a recent review).

Secondly, we contend that the problem we are dealing with should
really be regarded as one about (second-order) dynamical systems.  The
point is that a dynamical system may be represented as a system of
differential equations in many different coordinate
formulations; the question of real interest is whether there is {\em
some\/} representation of it which takes the form of an Euler-Lagrange
system with dissipation, not just whether a given representation of it
takes that form.  Of course this point applies equally, mutatis
mutandis, to the case in which there is no dissipation.  Now the
Helmholtz conditions as discussed in \cite{germ2,germ1}, in both the
classical and the generalized versions, suffer from the disadvantage
that they are conditions for a {\em given\/} system of differential
equations to be of Euler-Lagrange type.  There is, however, an
alternative approach to the problem which does deal with dynamical
systems rather than equations, at least in the case in which the
system can be expressed in normal form $\ddot{q}^a=F^a(\dot{q},q,t)$.
In this approach one asks (in the absence of dissipation) for
conditions for the existence of a so-called multiplier, a non-singular
matrix with elements $g_{ab}$, such that $g_{ab}(\ddot{q}^b-F^b)$
takes the Euler-Lagrange form (so that in particular when the
conditions are satisfied $g_{ab}$ will be the Hessian of the
Lagrangian with respect to the velocity variables).  The basic idea is
to put $h_a=-g_{ab}F^b$ in the conditions at the end of the
introduction, and regard the results as a system of partial
differential equations equations for $g_{ab}$ with $F^a$
known.  The seminal paper in this approach
is Douglas's of 1941 \cite{doug}, which analyses in great detail the
case of two degrees of freedom.  For a recent review of developments
since then see Sections 5 and 6 of \cite{olgageoff} and references therein.  One can in fact
also formulate conditions on a multiplier for a second-order dynamical
system, expressible in normal form, to be representable as equations
of Lagrangian form with dissipative forces of gradient type; these
generalize the known results for representation in Lagrangian form
without dissipation in an interesting way.

The new ingredient in \cite{germ2}, by comparison with \cite{germ1},
is the expression of the generalized Helmholtz conditions in terms of
quasi-velocities.  As presented in the paper this is quite a long-drawn-out
procedure, because in effect the conditions are rederived
from scratch.  Our third remark is that in principle this should be
unnecessary:\ a truly satisfactory formulation of the conditions
should be tensorial, in the sense of being independent of a choice of
coordinates (and of course quasi-velocities are just a certain type of
velocity coordinates).  The approach described in the previous
paragraph leads to conditions which have this desirable property.

Fourthly and finally, there is the question of whether generalized
Helmholtz conditions can be derived for other kinds of ``generalized
force'' terms than $\partial
D/\partial\dot{q}^a$.  One important case is that in which such a term
is of gyroscopic type.  We have obtained such conditions in this case,
again using the approach discussed in our second remark above.

These points are discussed in full detail in a recently written paper~\cite{MWC}.

\subsubsection*{Acknowledgements}
The first author is a Guest Professor at Ghent University:\ he is
grateful to the Department of Mathematical Physics and Astronomy at
Ghent for its hospitality.  The second author is a
Postdoctoral Fellow of the Research Foundation -- Flanders (FWO).

\end{document}